# SIMULATED POWER OF SOME DISCRETE GOODNESS-OF-FIT TEST STATISTICS FOR TESTING THE NULL HYPOTHESIS OF A 'ZIG-ZAG' DISTRIBUTION

## CLEMENT AMPADU\*, DANIEL WANG and MICHAEL STEELE[1]

Department of Mathematics

Central Michigan University

Mount Pleasant, MI, 48859, U. S. A.

e-mail: ampad1cb@cmich.edu

      wang1dx@cmich.edu

[1]Faculty of Business

Technology and Sustainable Development and

Faculty of Health Sciences and Medicine

Bond University

University Drive 4229, Robina, Queensland, Australia

e-mail: misteele@bond.edu.au

## Abstract

In this paper, we compare the powers of several discrete goodness-of-fit test statistics considered by Steele and Chaseling [10] under the null hypothesis of a 'zig-zag' distribution. The results suggest that the Discrete Kolmogorov-Smirnov test statistic is generally more powerful for the decreasing trend alternative. The Pearson Chi-Square statistic is generally more powerful for the increasing, unimodal, leptokurtic, platykurtic and bath-tub shaped alternatives. Finally, both the Nominal Kolmogorov-Smirnov and the Pearson Chi-Square test statistic are generally more powerful for the bimodal alternative. We also address the issue of the sensitivity of the test statistics to the alternatives under the 'zig-zag' null. In comparison to the uniform null of Steele and Chaseling [10], our investigation shows that the Discrete KS test statistic is most sensitive to the decreasing trend alternative; the Pearson Chi-Square statistic is most sensitive to both the leptokurtic and platykurtic trend alternatives. In particular, under the 'zig-zag' null we are able to clearly identify the most powerful test statistic for the platykurtic and leptokurtic alternatives, compared to the uniform null of Steele and Chaseling [10], which could not make such identification.

## 1. Introduction

The goodness-of-fit problem examines how well a sample of data agrees with a given distribution as its population. In the formal framework of hypothesis testing, the null hypothesis is that a given random variable follows a stated probability law, for example, the Beta-Normal distribution; the random variable may come from a process which is under investigation. The goodness-of-fit techniques applied to test the null hypothesis are based on measuring in some way the conformity of the sample data to the hypothesized distribution, or, equivalently, its discrepancy from it. The techniques usually give formal statistical tests and the measures of consistency or discrepancy are test statistics.

As a number of authors have noted including those in [10] little attention have been given in investigating the power of these test statistics, in particular those used in hypothesis tests of ordinal data. In [10] the authors compared the powers of six discrete goodness-of-fit test statistics for a uniform null distribution against a variety of fully specified alternative distributions. For the criteria used in their study, the results suggest that the test statistics specifically designed for ordinal data, namely, Discrete Kolmogorov-Smirnov (Discrete KS), Discrete Cramer-von Mises (Ordinal CVM), and Discrete Anderson Darling (Ordinal AD) are shown to be more powerful for the decreasing trend and step-type alternate distributions, whilst the nominal test statistics, namely, Pearson's Chi-Square and the Nominal Kolmogorov-Smirnov (Nominal KS), or the Watson (Ordinal Watson) test statistic are generally more powerful for the four other alternatives.

It is standard practice in the analysis of discrete power to hypothesize that the null distribution is uniform, since it is easy to explain the deviation in this case. In this paper, we take a different approach, and assume that the hypothesized distribution is "zig-zag," and we argue why shortly.

To our knowledge it appears that From [4] was the first to consider "zig-zag" as a trend alternative, moreover, a literature search also indicates that there does not appear to be a mathematical definition of "zig-zag", however, the reader will find the following definition useful.

**Definition 1.1.** A discrete probability distribution is said to be "*zig-zag*" if for $k$ categories the pattern of cell probabilities is one of the following:

(1) $p_1 < p_2 > p_3 < p_4 > \cdots < p_{k-1} > p_k$, provided that $k$ is odd,

(2) $p_1 > p_2 < p_3 > p_4 < \cdots > p_{k-1} < p_k$, provided that $k$ is odd,

(3) $p_1 < p_2 > p_3 < \cdots > p_{k-1} < p_k$, provided that $k$ is even,

(4) $p_1 > p_2 < p_3 > \cdots < p_{k-1} > p_k$, provided that $k$ is even,

where $0 \le p_i \le 1$ for $i = 1, ..., k$ and $\sum_{i=1}^{k} p_i = 1$.

An example of the zig-zag distribution is the null distribution used in this paper. Several more examples of the zig-zag distribution can also be found in From [4]. One can see that no formula is necessary to generate the "zig-zag" distribution. The requirements are minimal and they are (a) the pattern of cell probabilities must be one of the above (b) the cell probabilities must be chosen within the closed interval [0, 1] (c) all the cell probabilities must sum to one.

We now give our motivation for using zig-zag as the null distribution. Consider any type of data that can be characterized as ordinal. Let this data be a simple random sample of size $N$. Divide it up according to its $k$ distinct characteristics or categories. Let $n_i$ be the number of individuals in the sample with the $i$th characteristic, for $i = 1, ..., k$, where $\sum_{i=1}^{k} n_i = N$. The best estimate of the proportion of individuals with the $i$th characteristic is $\frac{n_i}{N}$ for $i = 1, ..., k$. Realistically, it is highly unlikely that the number of individuals in the sample with the $i$th characteristic is the same for $i = 1, ..., k$, that is, it is highly unlikely that $n_i \equiv n$ for $i = 1, ..., k$, and so $\dfrac{n}{N}$ is highly unlikely for $i = 1, ..., k$. The most realistic situation is that the proportion of individuals with the $i$th characteristic fluctuate for a majority of $i = 1, ..., k$. In light of this observation, $\dfrac{n}{N}$ is not a good estimate of the proportion of individuals with the $i$th characteristic for $i = 1, ..., k$.

Clearly, one can see from this argument that hypothesizing that the null distribution is zig-zag is a valid assumption, although a large power study will be needed to confirm this under various alternatives. However, as the reader will soon see, the mild result obtained in this paper, indicates that the use of the zig-zag null in favor of the uniform null is promising. Of course one could argue that it is possible to choose the cell probabilities so as to make the zig-zag distribution close to the uniform distribution, thus it would be desirable to use the uniform null rather than the zig-zag null. However, from a practical standpoint one can see that our argument is rigorous enough to accept the zig-zag null over the uniform null distribution. Even if a large power study refutes our argument, the result would still be subjective in nature, since simulated power is *not true* power.

Motivated by this argument, for the selected alternatives we examine the probability that the test statistics considered by Steele and Chaseling [10] will correctly lead to the rejection of a false null hypothesis under the zig-zag null.

Let $Z_i = \sum_{j=1}^{i} (O_j - E_j)$, $H_i = \sum_{j=1}^{i} E_j$, $\overline{Z} = \sum_{j=1}^{k} Z_j p_j$ and $N$ be the sample size. Further let $p_i$, $E_i$ and $O_i$ represent the cell probability, expected frequency, and observed frequency, in the $i$th cell, respectively. The test statistics considered by Steele and Chaseling [10] are in Table 1.

In Section 2, the alternative distributions are defined and the simulation and linear interpolation technique used to approximate the powers are discussed. In Section 3, we report the results of our power study for the zig-zag null against each alternative. We leave the reader with some concluding remarks in Section 4, in particular, we report on the sensitivity of the test statistics with respect to the zig-zag null of this paper and the uniform null considered by Steele and Chaseling [10].

**Table 1.** Test statistics used in the power study

| Name of test statistic | Equation | Author(s) |
|---|---|---|
| Pearson's Chi-square | $\chi^2 = \sum_{i=1}^{k} \dfrac{(O_i - E_i)^2}{E_i}$ | Pearson [6] |
| Discrete KS | $S = \max_{1 \le i \le k} \mid Z_i \mid$ | Pettitt and Stephens [7] |
| Ordinal CVM | $W^2 = N^{-1} \sum_{i=1}^{k} Z_i^2 p_i$ | Choulakian et al. [2] |
| Ordinal Watson | $U^2 = N^{-1} \sum_{i=1}^{k} (Z_i - \overline{Z})^2 p_i$ | Choulakian et al. [2] |
| Ordinal AD | $A^2 = N^{-1} \sum_{i=1}^{k} \dfrac{Z_i^2 p_i}{H_i(1 - H_i)}$ | Choulakian et al. [2] |
| Nominal KS | $NS = \dfrac{1}{2} \sum_{i=1}^{k} \mid O_i - E_i \mid$ | Pettitt and Stephens [7] |

**Definition 1.2.** By the sensitivity of a test statistic, we mean the probability of failing to detect the alternative trend, when such a trend is actually present. We shall define sensitivity by the formula $1 - p$, where $p$ is the simulated power of the test statistic.

**Remark 1.3.** It is evident from this definition that one could therefore consider sensitivity as *simulated type II error*, since the simulated power is not true power. It is also evident that chance of failing to detect the alternative trend, when such a trend is actually present, is low at higher simulated power levels, and high at lower simulated power levels. It is also clear that the most effective test statistic is the one with low sensitivity, that is, one whose simulated power is high.

## 2. The Power Study

Powers of each of the six test statistics are approximated for a zig-zag null across 10 cells against a selection of fully specified alternative distributions. Sample sizes used in this power study are 10, 20, 30, 50, 100 and 200, respectively.

The cell probabilities for the decreasing trend, bimodal distribution, leptokurtic distribution, and platykurtic distribution in this study are those used by Steele and Chaseling [10]. The cell probabilities for the increasing trend alternative were fully specified. The null distribution was also fully specified.

From [4] considered a new trend alternative derived from the Beta-Binomial distribution, this 'Beta-Binomial' trend alternative has cell probabilities:

$$P_i = \frac{\binom{k-1}{i-1}\Gamma(a+i-1)\Gamma(k+b-i)\Gamma(a+b)}{\Gamma(a+b+k-1)\Gamma(a)\Gamma(b)}, \quad i = 1, \ldots, k, \qquad (1)$$

where $\Gamma(\cdot)$ denotes the gamma function.

A discussion of the Beta-Binomial distribution can be found in Johnson et al. [5]. Denote equation (1) by the family $BB(a, b)$. For $a = b = 1$, we get the discrete uniform distribution. Of particular interest are the cases $a = b > 1$ which generate the unimodal symmetric alternative and the case $0 < a = b < 1$ which generates the bath-tub or U-shaped symmetric alternative. For the unimodal symmetric alternative, we fully specified the cell probabilities using $BB(1.5, 1.5)$ and for the bath-tub or U-shaped symmetric alternative we specified the cell probabilities using $BB(0.9, 0.9)$. Below is an example of the code in $R$ that was used to generate the ten cell probabilities for $BB(1.5, 1.5)$:

> $i = 1 : 10$

> $\dfrac{(\text{choose}(9, i-1) * \text{gamma}(0.5 + i) * \text{gamma}(11.5 - i) * \text{gamma}(3))}{(\text{gamma}(12) * \text{gamma}(1.5) * \text{gamma}(1.5))}$,

Damianou and Kemp [3], Steele and Chaseling [10] used the notion of linear interpolated power to derive the simulated power of these test statistics. A comparison of the relative power of categorical goodness-of-fit (GOF) test statistics is usually questionable if different significant levels are used for the null and alternative distributions of the test statistic. Let $\alpha$ be the desired level of significance, $\alpha_1$ be the significance level less than $\alpha$, and $\alpha_2$ be the significance level greater than $\alpha$, linear interpolation gives a weighting to the power based on how close $\alpha_1$ and $\alpha_2$ are to $\alpha$, the power is then computed using the formula:

$$\text{Power} = \frac{(\alpha - \alpha_1)P(T \geq X_2(\alpha) \mid H_1) + (\alpha_2 - \alpha)P(T \geq X_1(\alpha) \mid H_1)}{\alpha_2 - \alpha_1}, \qquad (2)$$

where $X_1(\alpha)$ and $X_2(\alpha)$ are the critical values immediately below and above the significance level $\alpha$, and $\alpha_1 = P(T \geq X_1(\alpha)|H_0)$ and $\alpha_2 = P(T \geq X_2(\alpha)|H_0)$ are the significance levels for $X_1(\alpha)$ and $X_2(\alpha)$, respectively.

For this study the power of each test statistic is estimated from 10,000 simulated random samples. The powers are obtained for critical values on both sides of the 1% level, and equation (2) is used to derive an approximate power for the 1% level. Whilst the 5% level is usually the norm, we selected the 1% level because it is more rigorous and tells us that there is only 1% chance that the simulated power results are due to chance, and thus the power by simulation is just as respectable as the power by calculation. The programming was done in the R language environment for statistical computing [8], and the code is available upon request from the first author. Table 2 summarizes the distributions used in our power study.

**Table 2.** Distributions used in the power study

| Description | Cell Probabilities | | | | | | | | | |
|---|---|---|---|---|---|---|---|---|---|---|
| | 1 | 2 | 3 | 4 | 5 | 6 | 7 | 8 | 9 | 10 |
| Zig-Zag (Null Distribution) | 0.20 | 0.05 | 0.10 | 0.05 | 0.10 | 0.02 | 0.20 | 0.10 | 0.08 | 0.10 |
| Decreasing | 0.32 | 0.13 | 0.10 | 0.08 | 0.07 | 0.07 | 0.06 | 0.06 | 0.05 | 0.05 |
| Increasing | 0.03 | 0.04 | 0.05 | 0.06 | 0.10 | 0.11 | 0.12 | 0.14 | 0.16 | 0.19 |
| Unimodal | 0.06 | 0.09 | 0.17 | 0.17 | 0.12 | 0.12 | 0.12 | 0.17 | 0.09 | 0.06 |
| Bimodal | 0.05 | 0.11 | 0.17 | 0.11 | 0.06 | 0.06 | 0.11 | 0.17 | 0.11 | 0.05 |
| Leptokurtic | 0.05 | 0.05 | 0.05 | 0.05 | 0.30 | 0.30 | 0.05 | 0.05 | 0.05 | 0.05 |
| Platykurtic | 0.04 | 0.11 | 0.11 | 0.12 | 0.12 | 0.12 | 0.12 | 0.11 | 0.11 | 0.04 |
| Bath-tub Shaped | 0.11 | 0.10 | 0.10 | 0.01 | 0.09 | 0.09 | 0.10 | 0.10 | 0.10 | 0.11 |

### 3. Powers of the Test Statistic for Each Alternative[1]

### 3.1. Decreasing trend alternative distribution

With the decreasing trend alternative, Figure 1 indicates that the Discrete KS statistic has the highest power, followed by the Pearson Chi-Square. When the sample size is at most 30, the Ordinal CVM statistic and the Nominal KS test statistic have comparable power. When the sample size is at least 30, the Nominal KS has higher power compared to the Ordinal CVM test statistic. The Ordinal Watson test statistic has superior power compared to the Ordinal AD test statistic.

---

[1]It should be noted that the power analysis for Section 3 was done from power charts produced in Microsoft Excel 2003. However, for publication purposes the graphs have been formatted using MINITAB 15.

The Ordinal AD test statistic has the poorest power in comparison to the powers of the other test statistics under this alternative. When the sample size is at least 100, all the statistics have power very close to 1.

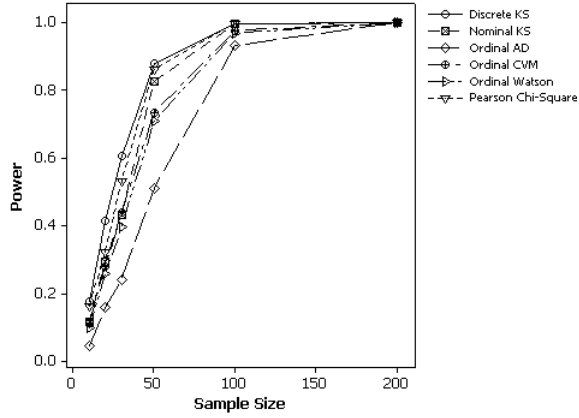

**Figure 1.** Power of the test statistics for a zig-zag null and decreasing alternative.

### 3.2. Increasing trend alternative distribution

With the increasing trend alternative, Figure 2 indicates that the Pearson Chi-Square statistic has the highest power. The powers of the Ordinal Watson test statistic and the Ordinal CVM test statistic are similar. The powers of the Nominal KS test statistic dominate those of the Discrete KS statistic under this alternative. The Ordinal AD test statistic has the poorest power in comparison to all the test statistics under this alternative. When the sample size is at least 100, all the statistics have power very close to 1.

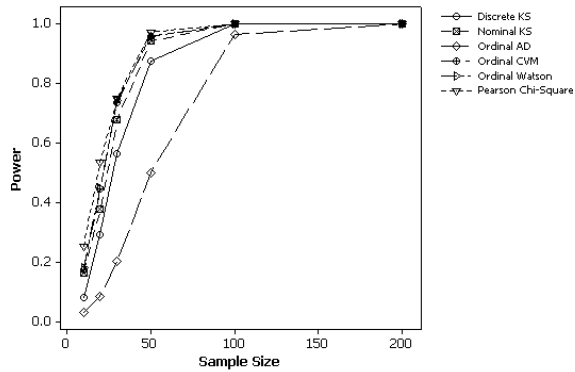

**Figure 2.** Power of the test statistics for a zig-zag null and increasing alternative.

### 3.3. Unimodal alternative distribution

With the unimodal alternative, Figure 3 indicates that the Pearson Chi-Square statistic has the highest power followed by the Nominal KS, Ordinal AD, and the Discrete KS test statistic, respectively. When the sample size is at most 30, the Ordinal Watson and the Ordinal CVM statistic appear to have similar power, and when the sample size is at least 30 the Ordinal Watson test statistic appear to have superior power in comparison to the Ordinal CVM statistic. When the sample size is at most 50, the powers of the Discrete KS, Ordinal Watson, and the Ordinal CVM statistic are not meaningful. When the sample size is at least 100, the Pearson Chi-Square, Nominal KS and the Ordinal AD test statistics have power very close to 1.

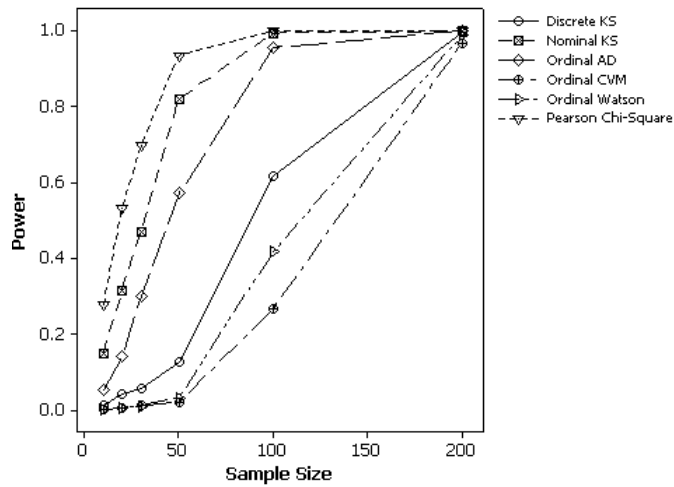

**Figure 3.** Power of the test statistics for a zig-zag null and unimodal alternative.

### 3.4. Bimodal alternative distribution

With the bimodal alternative, Figure 4 indicates that the Pearson Chi-Square and the Nominal KS statistic have comparable power and are the highest under this alternative, as one can see there is also a significant difference when their powers are compared with the others. When the sample size is at most 30, the powers of the Ordinal Watson, Discrete KS, and the Ordinal CVM statistics are approximately the same. When the sample size is at least 30 but no more than 50, the Ordinal Watson and Discrete KS statistic have similar and higher power in comparison to the Ordinal CVM test statistic. When the sample size is at least 50, the Discrete KS test statistic has greater power compared to the Ordinal Watson and Ordinal CVM test statistic.

When the restriction on the sample size is at most 30, the Ordinal Watson, Discrete KS and the Ordinal CVM test statistic have the worst power compared to the other test statistics. On the other hand when the restriction on the sample size is at least 30, the Ordinal CVM statistic has the worst power. When the sample size is at least 100, the Pearson Chi-Square, Nominal KS and the Ordinal AD test statistics have power very close to 1.

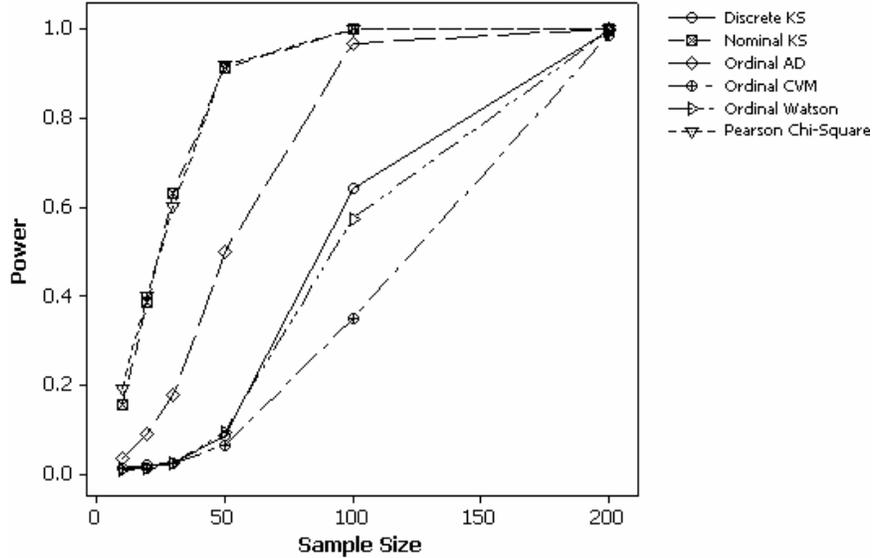

**Figure 4.** Power of the test statistics for a zig-zag null and bimodal alternative.

### 3.5. Leptokurtic alternative distribution

With the leptokurtic trend alternative, Figure 5 shows that the Pearson Chi-Square test statistic has the highest power when the sample size is no more than 20. When the sample size is at least 20, the powers of the Pearson Chi-Square and the Nominal KS test statistic converge to 1 quite rapidly and are the highest compared to the other test statistics. The Ordinal AD test statistic has higher power compared to the Discrete KS statistic when the sample size is at most 50, otherwise we see that the powers of these test statistics converge to 1. The Ordinal Watson and Ordinal CVM test statistic have similar power across the sample size, and their powers are the worst in comparison to the other test statistic under this alternative. It also appears that the powers of all the test statistics are very close or exactly 1 when the sample size is at least 100.

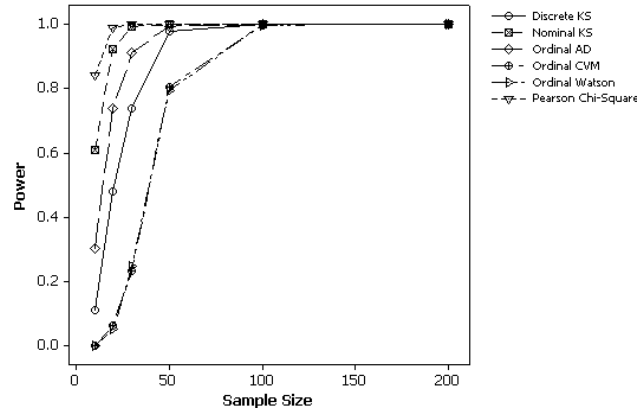

**Figure 5.** Power of the test statistics for a zig-zag null and leptokurtic alternative.

### 3.6. Platykurtic alternative distribution

With the platykurtic trend alternative, Figure 6 indicates that the Pearson Chi-Square statistic has the highest power followed by the Nominal KS, Ordinal AD, and the Discrete KS test statistic, respectively. When the sample size is at most 30, the Ordinal Watson and Ordinal CVM test statistic have approximately the same power, and their powers are the worst compared to the other test statistics. On the other hand when the sample size is at least 30, the Ordinal Watson test statistic has superior power compared to the Ordinal CVM test statistic. The Ordinal CVM test statistic has the worst power when the sample size is at least 30. When the sample size is at least 100, the powers of Pearson Chi-Square, Nominal KS and Ordinal AD tend to 1.

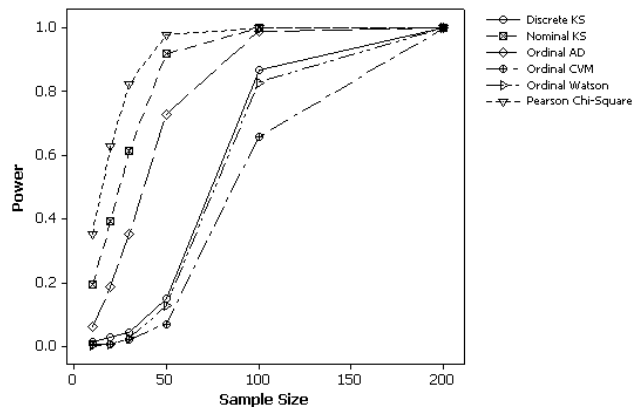

**Figure 6.** Power of the test statistics for a zig-zag null and platykurtic alternative.

### 3.7. Bath-tub shaped alternative distribution

With the bath-tub shaped trend alternative, Figure 7 indicates Pearson Chi-Square test statistic has the highest power followed by the Nominal KS statistic. When the sample size is at most 20, the Discrete KS, Ordinal CVM, Ordinal AD, and Ordinal Watson test statistics have approximately the same power and are worst in comparison to the Pearson Chi-Square and Nominal KS test statistics. When the sample size is at least 20, the Ordinal AD has the highest power compared to the Discrete KS, Ordinal Watson and the Ordinal CVM test statistics, respectively. Also when the sample size is at least 20, the Ordinal CVM has the worst power compared to all the other statistics. It also seems that for the Ordinal Watson, Ordinal CVM, Ordinal AD and the Discrete KS statistic, the powers never surpass 0.55 even when the sample size is very large. We also find that when the sample size is least 100, the Pearson Chi-Square and the Nominal KS statistic have powers tending to 1.

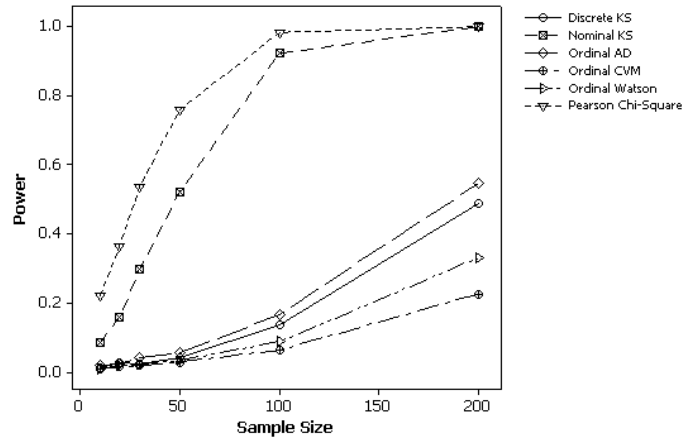

**Figure 7.** Power of the test statistics for a zig-zag null and bath-tub shaped alternative.

## 4. Conclusions

For testing the hypothesis of a zig-zag null, we find the following: (1) The Discrete KS test statistic has the highest power for the decreasing trend alternative. (2) The Pearson Chi-Square test statistic has the highest power for the increasing, unimodal, leptokurtic, platykurtic and bath-tub shaped alternatives. (3) The Pearson Chi-Square test statistic and the Nominal KS test statistic have comparable power and are the highest for the bimodal trend alternative.

Further, although some of the alternatives considered in this paper were not considered by Steele and Chaseling [10] and our simulated powers were obtained at the 1% significance level rather than the 5% level as they did, we find for those alternatives considered in both papers the following: (1) When the trend is decreasing, and the null hypothesis is uniform, Steele and Chaseling [10] show that the power of the Nominal KS statistic lags behind the rest. However, in the present paper, under the zig-zag null, the Nominal KS statistic has approximately the third highest power. (2) When the trend is platykurtic, and the null hypothesis is uniform, the power results of Steele and Chaseling [10] show that there is a conflict as to which test statistic has optimum power and can therefore detect the trend correctly. However, in the present paper, under the zig-zag null, the Pearson Chi-Square statistic wins out. Also under the uniform null, Steele and Chaseling [10] show that the power of the Discrete KS statistic lags behind the rest, however, in the present paper, under the zig-zag null the same statistic has the fourth highest power. Also, under the uniform null when the sample size is at most 3 per cell, the power results of Steele and Chaseling [10] are not meaningful for the Ordinal AD, Discrete KS and the Discrete CVM test statistics. However, under the zig-zag null when the sample size is at most 30, one clearly sees that this trend has been reversed at least for the Ordinal AD test statistic.

Finally, under the uniform null of Steele and Chaseling [10] when the sample size is very large, it appears that the Ordinal AD, Discrete KS and the Discrete CVM test statistic have powers that do not surpass 0.60, however, under the zig-zag null, these statistics have power tending to 1. (3) When the trend is leptokurtic, and the null hypothesis is uniform, the power study of Steele and Chaseling [10] show there is a conflict as to which test statistic has optimum power and can detect the trend correctly. However, in the present paper, under the zig-zag null the Pearson Chi-Square statistic wins out. Also from the power study of Steele and Chaseling [10], under the uniform null the Ordinal AD test statistic lags behind the rest. However, in the present paper, under the zig-zag null the Ordinal AD test statistic has the third highest power. (4) When the trend is Bimodal, and the null hypothesis is uniform, the power results of Steele and Chaseling [10] show that the Ordinal AD, Discrete KS, and the Discrete CVM test statistic give power that do not surpass 0.50 even for the large sample sizes considered, however in the present paper, under the zig-zag null, the powers of these statistics tend to 1.

Table 3 summarizes the main differences in the sensitivity of the test statistics under the uniform null of Steele and Chaseling [10] and the zig-zag null of the present paper. For example, when the trend is platykurtic and the null hypothesis is uniform the most sensitive statistic could not be identified, in particular the power results of Steele and Chaseling [10] show that Discrete Watson, Pearson Chi-Square and Nominal KS have approximately the same power. However under the zig-zag null, the Pearson Chi-Square is most sensitive by Definition 1.2 and Remark 1.3.

**Table 3.** Sensitivity comparison under the uniform and zig-zag null

| Alternative | Uniform Null | Zig-Zag Null |
| --- | --- | --- |
| Decreasing | Discrete KS is the third most sensitive | Discrete KS is most sensitive |
| Platykurtic | Most sensitive statistic cannot be determined | Pearson Chi-Square is most sensitive |
| Leptokurtic | Most sensitive statistic cannot be determined | Pearson Chi-Square is most sensitive |

From this study there is evidence under the zig-zag null that some of the test statistics are sensitive to the alternatives considered in both papers. It appears sensitivity analysis of categorical goodness-of-fit test statistics to various alternatives needs further investigation, perhaps in the spirit of Basu et al. [1]. Basu et al. [1] investigated the sensitivity of the power of the family of divergence test statistics constructed by Read [9] for the "bump" and "dip" alternatives. Basu et al. [1] note an improvement in the power of the family of these statistics when the issue of sensitivity is considered in the construction of some new test statistics which also belong to the family of divergence test statistics of Read [9].